\documentclass{elsart}

\usepackage{amssymb}
\usepackage{amsmath}
\usepackage{amsfonts}
\usepackage{diagrams}

\def\aa{{\mathcal A}}
\def\KK{{\mathbb K}}
\def\AA{{\mathcal A}}
\def\ii{{\mathcal I}}
\def\oo{{\mathcal O}}
\def\MM{{\mathbb M}}

\def\CC{{\mathbb C}}
\def\PP{{\mathbb P}}
\def\NN{{\mathbb N}}
\def\QQ{{\mathbb Q}}
\def\RR{{\mathbb R}}
\def\ZZ{{\mathbb Z}}
\def\QEI{{Q\setminus E_I}}

\begin{document}

\begin{frontmatter}

% Title, authors and addresses

% use the thanksref command within \title, \author or \address for footnotes;
% use the corauthref command within \author for corresponding author footnotes;
% use the ead command for the email address,
% and the form \ead[url] for the home page:

\title{Implicitization of rational surfaces using toric
varieties\thanksref{label1}} \thanks[label1]{Amit Khetan was supported
by an NSF postdoctoral fellowship (DMS-0303292). Carlos D'Andrea was
supported by a Miller Research Fellowship (2002-2005).}

\author{Amit Khetan} 
\ead{khetan@math.umass.edu}
\ead[url]{http://www.math.umass.edu/~khetan}
% \corauth[cor1]{}
\address{Department of Mathematics \\ University of Massachusetts at Amherst\\  Amherst MA 01002, USA}

\author{Carlos D'Andrea}
\ead{cdandrea@math.berkeley.edu}
\ead[url]{http://www.math.berkeley.edu/~cdandrea}
\address{Miller Institute for Basic Research in Science and Department of Mathematics \\ University of California at Berkeley \\ Berkeley, CA 94720, USA}

\begin{abstract}
A parameterized surface can be represented as a projection from a
certain toric surface. This generalizes the classical homogeneous and
bihomogeneous parameterizations. We extend to the toric case two
methods for computing the implicit equation of such a rational
parameterized surface. The first approach uses resultant matrices and
gives an exact determinantal formula for the implicit equation if the
parameterization has no base points. In the case the base points are
isolated local complete intersections, we show that the implicit
equation can still be recovered by computing any non-zero maximal
minor of this matrix.

The second method is the toric extension of the method of moving
surfaces, and involves finding linear and quadratic relations
(syzygies) among the input polynomials. When there are no base
points, we show that these can be put together into a square
matrix whose determinant is the implicit equation. Its extension
to the case where there are base points is also explored.
\end{abstract}

\end{frontmatter}

% main text
\section{Introduction}
A rationally parameterized surface $\Phi(s,t)$ in affine three
space is defined by a map $\phi : \CC^2 \to \CC^3$ given by three
rational components:

\begin{equation} \label{par}
\left\{
\begin{array}{ccc}
X_1&=&\frac{x_1(s,t)}{x_4(s,t)} \\
X_2&=&\frac{x_2(s,t)}{x_4(s,t)} \\
X_3&=&\frac{x_3(s,t)}{x_4(s,t)}
\end{array}
\right.
\end{equation}

Here $x_1,x_2,x_3,x_4$ are (Laurent) polynomials in two variables $s$
and $t$ with coefficients in $\CC$ (or $\RR$ or $\QQ$ or any arbitrary
subfield $\KK$ of $\CC$). Let $\Phi \subset \CC^3$ be the
smallest algebraic surface containing (\ref{par}). The
\textit{implicitization problem} \cite{CLO1,cox1} is to compute the
polynomial equation $P(X_1, X_2, X_3)$ defining $\Phi$.  At times we
will also consider the same surface in projective space, where there
are four coordinates with equations given by $X_1, X_2, X_3,$ and
$X_4$.  \par

The last few decades have witnessed a rise of interest in the
implicitization problem for geometric objects motivated by
applications in computer aided geometric design and geometric
modelling (\cite{AGR,AS,BCD,buc,CGKW,cox1,CGZ,dok,kal,MC,SAG}).  A
very common approach is to write the implicit equation as the
determinant of a matrix whose entries are easy to compute.  \par

Our approach is also to look for matrix formulas, but we recast the
parameterization in terms of a projection from a certain toric
surface built out of the specific monomials which appear in $x_1,
x_2, x_3, x_4$. This generalizes the standard approaches of
projections from tensor product surfaces (Segre embeddings of
$\PP^1 \times \PP^1$) or from total degree surfaces (Veronese
embeddings of $\PP^2$).  So while previously $x_1, x_2, x_3, x_4$
have been considered only as ``generic'' homogeneous or
bihomogeneous polynomials, we can exploit sparsity present in the
parameterization.\par

Standard homogenization of sparse polynomials can result in numerous
spurious base points of the projection at infinity. By using the
more general toric surface, customized for the equations on hand, many
of these extraneous base points at infinity can be avoided. This
results in smaller matrices and fewer extraneous factors in the
computation of the implicit equation. Toric projections can also be
exploited in the construction of the parameterization. The work of
Krasauskas \cite{Kra} shows how ``toric surface patches'' can be used
to parametrize regions on a surface shaped like arbitrary sided
polygons. 

In this article we extend to the toric case two methods for
computing the implicit equation:  computing a Chow form and
computing syzygies on the input polynomials $x_1, x_2, x_3, x_4$.

A classical method for finding the implicit equation is to compute the
bivariate resultant or Chow form of the three polynomials

\begin{eqnarray} \label{pols}
f_1 &=& x_1(s,t)-X_1\,x_4(s,t) \notag\\
f_2 &=& x_2(s,t)-X_2\,x_4(s,t) \notag\\
f_3 &=& x_3(s,t)-X_3\,x_4(s,t) 
\end{eqnarray}

Our first approach essentially follows the classical method using
the sparse resultant in place of the classical bivariate
resultant. Formally, we will reduce the computation of $F$ to the
computation of the Chow form of a toric surface which projects
onto $\Phi.$ Exact matrix formulas for computing this Chow form
were found by the first author in \cite{khe}.

\par We show that if the projection has no {\em base points}, points on the
toric variety such that $f_1, f_2, f_3$ are simultaneously zero, the matrix
constructed gives an exact determinantal formula for the implicit
equation. New to our approach is an analysis when base points are
present. We show that if the base points are isolated local complete
intersections, the implicit equation can still be recovered by
computing a non-zero maximal minor of this matrix.

The second method involves finding linear and quadratic relations
(syzygies) among the polynomials $x_1,x_2,x_3,x_4$ of a certain fixed
type.  When there are no base points, we will see how these can be put
together into a square matrix whose determinant is exactly the
implicit equation. This is precisely the technique used in the method
of moving surfaces for tensor product or total degree surfaces
(\cite{SC,BCD,CGZ,cox1,SAG}). Our contribution is to exploit the
structure of the sparsity of the polynomials to avoid extra base
points.  Moreover, we present a novel proof of the validity of the
method of moving surfaces which ties together the complexes of moving
planes and quadrics with the resultant complex in a natural way. This is
described in Section 5.

The method of moving surfaces can also be applied in the presence of
basepoints and often still produces the correct implicit
equation. Validity in the presence of basepoints was proved under
certain conditions in the total degree \cite{BCD} and tensor product
\cite{AHW} situations. We do not have a proof in the general toric
setting but we illustrate the situation with a few examples. 

The paper is organized as follows: in Section \ref{s:2}, we recall
some properties of toric surfaces and introduce some notation.  In
Section \ref{s:3}, we present the first of our methods and show
that it works if the base points are a local complete
intersection. Next, we present in Section \ref{s:4} the method of
moving quadrics and shows its validity in the case where there are
no base points. We also give some examples and an exploration of
what happens when base points are present.

\section{Toric varieties, parameterizations, and base points} \label{s:2}

Let $\aa = \{ \alpha_1, \dots, \alpha_N\} \subset\ZZ^2,$ a finite
subset of points, and $Q$ the convex hull of the points in $\aa$. The
{\em toric variety} $X_{\aa}$ associated with $\aa$ is defined as the
Zariski closure of the set of points $(x^{\alpha_1}: \dots:
x^{\alpha_N})$ in $\PP^{N-1}$ where $x$ ranges over $(\CC^\ast)^2$
(the ``algebraic'' torus).  See \cite{GKZ,CLO} for details.

If each of the polynomials $x_i$ has its support contained in $\aa,$
then it is a linear combination of monomials in $\aa,$ hence can be
thought of as a linear functional on $\PP^{N-1}$ defining a hyperplane
section of $X_{\aa}$.  If each $x_i$ has a different support $\aa_i$
we will define $\aa$ as the union of the supports; that is:

$$ \aa = \aa_1 \cup \aa_2 \cup \aa_3 \cup \aa_4 $$

Therefore we can consider the set of zeros $Z$ of $x_1, x_2, x_3, x_4$
in $X_{\aa}$.  Note that $Z$ will contain those common zeros of the
$x_i$ in $(\CC^{\ast})^2$.  Now the map $\phi$ can be realized as (the
affine part of) a projection from $X_{\aa}$ to $\PP^3$ via the
hyperplane sections $x_1, x_2, x_3, x_4$.

The points in $Z$ correspond to basepoints of this projection. We
will assume that $\gcd(x_1,x_2,x_3,x_4)=1,$ so that $Z$ is finite.
For each $p\in Z,$ we get a certain multiplicity
$e(\ii_{Z,p},\oo_{X_\aa,p}).$ The degree of the parameterization
$\phi$ is the generic number of points in $X_{\aa}$ which map to a
point in $\Phi$. The degree of $\Phi$ is the total degree of its
implicit equation.

\par Now, as in \cite[Appendix]{cox1}, we have the following
\textit{degree formula:}

\begin{prop}\label{df}
$$\deg(\phi)\deg(\Phi)={\rm Area}(Q)-\sum_{p\in Z}e(\ii_{Z,p},\oo_{X_\aa,p})$$
where ${\rm Area}(Q)$ is the normalized area of the polygon $Q$
equal to twice its usual Euclidean area (in particular ${\rm
Area(Q)}$ is always an integer).
\end{prop}

In the next section we will consider the case when there are no
basepoints, that is $Z = \emptyset$. If $x_1, x_2, x_3, x_4$ are each
generic with respect to their supports $\aa_1, \aa_2, \aa_3, \aa_4$
then this will be the case provided a certain geometric condition
on the supports $\aa_i$ holds. This is expressed in the next
result.

\begin{prop} \label{genericnobp}
Fix subsets $\aa_1, \aa_2, \dots, \aa_k \subset \ZZ^2$ with $k \geq 3$
and define $\aa = \aa_1 \cup \cdots \cup \aa_k$. Let $x_1, \dots, x_k$
be generic polynomials with $x_i$ supported on $\aa_i$.  Let $Q_i =
Conv(\aa_i)$ and $Q = Conv(\aa)$ be the associated polytopes. Assume
that $\dim(Q) = 2$. The polynomials $x_i$ viewed as sections of
$X_{\aa}$ have no common zeros if and only if every edge of $Q$
intersects at least two of the polytopes $Q_i$.
\end{prop}

\begin{pf} The torus orbits of $X_{\aa}$ correspond to the faces
of $Q_{\aa}$. The restriction of $x_i$ to a particular orbit
corresponds to intersecting $Q_i$ with the corresponding face of $Q$
and setting all terms of $x_i$ not in the intersection to 0.  A
zero-dimensional face of $Q$ corresponds to a vertex, which by
construction must be a vertex of some $Q_i$. The corresponding $x_i$
restricts to a single non-zero monomial with a generic (non-zero)
coefficient, hence does not vanish at this point. On the orbit
corresponding to the dense 2 dimensional torus, any 3 of the
polynomials do not generically have a common zero. Finally, a
one-dimensional orbit corresponds to an edge of $Q$. By hypothesis,
intersecting with the $Q_i$ yields at least two non-zero polynomials
with generic coefficients which do not have a common zero on the
one-dimensional space. Conversely, if an edge of $Q$ intersects only
one of the $Q_i$ then it must be an edge of that $Q_i$, so that $x_i$
restricts to a polynomial in one variable, while all other $x_j$
restrict to zero. Thus every root of this restriction of $x_i$ is a
common zero of all of the $x_j$.

\end{pf}

Of course if all of the $x_i$ have the same support, the case most
often of interest for implicitization, then the condition above is
automatically satisfied. In general it corresponds to a mild geometric
compatibility of the supports.

\section{Implicitization from the Chow Form}\label{s:3}

The {\em Chow form} of $X_{\aa}$ is a polynomial $Ch_\aa$ in the
coefficients of three linear sections $f_1, f_2, f_3$ which is zero
whenever $f_1, f_2, f_3$ have a common root on $X_{\aa}$. In the case
where $Z=\emptyset,$ we get the following result.

\begin{thm} \label{chow}
Let $x_1,x_2,x_3,x_4$ be Laurent polynomials with complex
coefficients.  Let $Z$ be the set of common zeros on the toric variety
$X_{\aa}$ corresponding to the union of their supports. Let $Ch_\aa$
be the Chow form of the toric variety $X_\aa.$ Let $f_1, f_2, f_3$ be
as in (\ref{pols}) and $P(X_1, X_2, X_3)$ the implicit equation of
$\Phi.$ If $Z = \emptyset$, then there exists a nonzero constant
$c\in\CC$ such that

\begin{equation}\label{cuco}
Ch_\aa(f_1, f_2, f_3)=cP^{\deg(\phi)}.
\end{equation}

\end{thm}

\begin{pf}
Let $G(X_1,X_2,X_3)$ be the left-hand side of \eqref{cuco}. For a
generic point on the surface there is an associated common zero of
$f_1, f_2, f_3$. Conversely if $X_1, X_2, X_3$ are such that $f_1,
f_2, f_3$ have a common zero $(s,t)$ then as $Z = \emptyset$,
$x_4(s,t) \neq 0$ and thus $(X_1, X_2, X_3)$ is a point on the
surface. As $P$ is irreducible, it follows that $G=cP^d,$ with $c
\neq 0$ and $d\in\NN.$ In order to verify that $d=\deg(\phi),$ it
is enough to see that the degree of $P$ is $\rm{Area}(Q).$ This
follows easily by applying the Chow form to the dual Pl\"ucker
coordinates of the polynomials (\ref{pols}) (see \cite{GKZ}) and
by noting that the dual Pl\"ucker coordinates have degree one in
$X_1,X_2,X_3.$ The degree of $Ch_\aa$ in the Pl\"ucker coordinates
is ${\rm Area}(Q).$
\end{pf}

In \cite{khe2} there is a construction for computing the Chow form
of any toric surface. Given a toric surface $X_\aa$ with $Q = {\rm
conv}(\aa)$, and three sections $f_1, f_2, f_3$ with $f_i =
\sum_{a \in \aa} C_{ia}x^{\alpha}$. Then $Ch_{\aa}(f_1,
f_2, f_3)$ is the determinant of a matrix of the following block
form:

\begin{equation*}
\label{e:blmtrx}
\begin{pmatrix}
 B & L \\ \tilde{L} & 0 \end{pmatrix},
\end{equation*}
\noindent

Here the entries of $L$ and $\tilde{L}$ are linear forms, and the
entries of $B$ are cubic forms in the coefficients $C_{ia}$, as
described below.

The columns of $B$ and $\tilde{L}$ are indexed by the lattice
points in $Q$, the rows of $B$ and $L$ are indexed by the interior
lattice points in $2 \cdot Q$, the matrix $\tilde{L}$ has three
rows indexed by $\{ f_1, f_2, f_3 \}$, and the columns of the
matrix $L$ are indexed by pairs $(f_i, a)$ where $i \in \{1,2,3\}$
and $a$ runs over the interior lattice points of $Q$.  Each entry
of $L$ and $\tilde L$ is either zero or is a coefficient of some
$f_i$ and is determined in the following straightforward manner.
The entry of $\tilde L$ in row $f_i$ and column $a$ is the
coefficient of $x^a$ in $f_i$. The entry of $L$ in row $b$ and
column $(f_i, a)$ is the coefficient of $x^{b-a}$ in $f_i$. The
entries of the matrix $B$ are linear forms in {\em bracket
variables}.  A bracket variable is defined as

$$[abc] = \det \bmatrix  C_{1a} & C_{1b} & C_{1c} \\
                         C_{2a} & C_{2b} & C_{2c} \\
                         C_{3a} & C_{3b} & C_{3c} \endbmatrix,$$

\noindent There is an explicit, combinatorial construction of the
matrix $B$ given in \cite{khe}. By virtue of Theorem \ref{chow}
above we get the immediate corollary.

\begin{cor}
If $Z=\emptyset,$ then there is a determinantal formula $M_{\aa}$
for computing $P^{\deg(\phi)}.$
\end{cor}

\begin{exmp} \label{ex1}

 Consider the surface parameterized by

\begin{eqnarray*}
x_1 &=& s^3 + t^2\\
x_2 &=& s^2 + t^3\\
x_3 &=& s^2t + st^2\\
x_4 &=& st
\end{eqnarray*}

The associated polygon $Q$ is a quadrilateral in the first quadrant
with vertices $(2,0), (3,0), (0,2), (0,3)$. So we compute the Chow
form, with respect to this polygon of $s^3 + t^2 - X_1st, s^2 + t^3 -
X_2st, s^2t + st^2 - X_3st$ which results in the following $7 \times
7$ matrix:

$$
\begin{bmatrix}
0 & 1 & -X_1  & 0     & 1     & 0    & 0 \\
1 & 0 & -X_2  & 0     & 0     & 0    & 1 \\
0 & 0 & -X_3  & 1     & 0     & 1    & 0 \\
0 & 0 & 1     & -X_3  & -X_2  & X_1  & X_1 - 1 \\
0 & 0 & 0     & 0     & X_3   & -1   & -1  \\
0 & 0 & -X_1  & X_2-1 & 1     & X_2  & -X_3 \\
0 & 0 & 1-X_1 & X_2   & 1-X_2 & -X_3 & X_1
\end{bmatrix}
$$

The determinant of this matrix is:

\begin{eqnarray*}
&&2 +X_1-5X_3^3-X_1^2X_2-X_2^2X_1+X_3X_1^3+X_3X_2^3+X_3^5+X_2+5X_3\\
&&+4X_2^2X_3^2-X_3X_1-2X_2X_1-X_2X_3-3X_2X_3^2+X_2^2X_3-3X_2X_1X_3^3\\
&&-X_2^2X_3^2X_1+4X_1^2X_3^2 -3X_1X_3^2+X_2^2X_1X_3+X_3X_1^2X_2-X_1^2X_3^2X_2\\
&&+2X_2X_3^3 + X_1^2X_3 - 5X_1X_2 X_3 + 2X_1X_3^2
\end{eqnarray*}

This is the degree $5$ (equal to ${\rm Area}(Q)$) affine implicit
equation.
\end{exmp}

It is an immediate consequence from the proof of Theorem
\ref{chow} that when $Z \neq \emptyset$ the Chow form $Ch_\aa$ is
identically zero. However, we shall see in the next section that
the implicit equation can still be recovered from maximal minors
of the resultant matrix. This shows how a matrix resultant formula
encodes much more information than just the Chow form.

\subsection{Base Points}

In this section we take a closer look at the Chow form matrix
described above in order to determine what happens in the presence of
base points. Throughout this section we will assume $x_1, x_2, x_3,
x_4$ are specific choices of polynomials supported on $\aa$ which may
in particular have base points.

We will see that we can always get a matrix whose determinant is a
non-trivial multiple of the implicit equation.  In order to still get
an exact formula we will need a hypothesis on the structure of the
basepoints.  By the construction of $\aa$, we will always be able to
assume that the points in $Z$ are smooth points of $X$ (see the explanation
in the proof of Theorem \ref{impbp}). In that case
the local ring $\oo_{X,p}$ is just the localized polynomial ring in
two variables $x, y$.

\begin{defn} Let $X$ be a variety of dimension $n$. A {\em zero
dimensional local complete intersection (LCI)} is a subscheme $Z$ in
the smooth locus of $X$, such that for each point $p$ in $Z$, the
ideal $I_{Z,p}$ of the local ring $\oo_{X,p}$ is defined by $n$
equations.
\end{defn}

The main property of local complete intersections that we will use is
contained in the next proposition.

\begin{prop}
If $Z$ is a local complete intersection then the multiplicity \\
$e(\ii_{Z,p},\oo_{X,p})$ is equal to the vector space dimension of the
finite local algebra $\oo_{Z,p} = \oo_{X,p}/I_{Z,p}$. In particular
$\sum_{p \in Z} e(\ii_{Z,p},\oo_{X,p})$ is equal to the vector space
dimension of the affine coordinate ring of $Z$.
\end{prop}

This proposition is a consequence of \cite[Theorem 4.7.4]{BH} as
$I_{Z,p}$ is generated by a regular sequence. Hence, the Euler
characteristic is just the length of $\oo_{X,p}/I_{Z,p}$, which
is the vector space dimension in the zero-dimensional case.

We can now state the main result of this section:

\begin{thm}\label{impbp}

Let $\pi : X_{\aa} \to \PP^3$ be a projection onto a surface $\Phi$
parameterized by $x_1, x_2, x_3, x_4$ with no common factor such that
$\aa$ is the union of the supports of the $x_i$. Let $Z \subset
X_{\aa}$ be the finite set of basepoints of $\pi$. Now, let $M_{\aa}$
be the determinantal formula for $Ch_{\aa}$ from \cite{khe}. where
$f_1, f_2, f_3$ are the polynomials $x_1(s,t) - X_1x_4(s,t), x_2(s,t)
- X_2x_4(s,t),$ and $x_3(s,t) - X_3x_4(s,t)$ respectively. Then the
implicit equation $P^{\deg(\phi)}$ divides any maximal minor of
$M_{\aa}$.

Moreover,a maximal minor of $Ch_{\aa}(f_1, f_2, f_3)$ using all of
the Sylvester rows and columns exists and has determinant equal to
exactly $P^{\deg(\phi)}$ if:

\begin{enumerate}
\item $Z$ is a local complete intersection on $X_{\aa}$.
\item The Sylvester columns in $L$, indexed by $\rm{int}(Q)$, are
linearly independent for generic choices of $X_1, X_2,
X_3$. Equivalently, $f_1, f_2, f_3$ have no syzygies supported on
$\rm{int}(Q)$ with coefficients in $\CC[X_1,X_2,X_3]$.
\end{enumerate}
\end{thm}

The LCI hypothesis seems to be ubiquitous in implicitization \cite{BC,BCD}.
Note that in particular isolated basepoints are always LCI so
that ``generically'' even if $x_1, x_2, x_3, x_4$ have basepoints,
e.g. if the geometric condition of Theorem \ref{genericnobp} is not
satisfied, the basepoints they do have will be LCI.

The second condition is somewhat more subtle and is not really well
understood except that it was quite difficult to construct examples
for which it fails (see Example 3.5). It can be compared with the
Assumption 5.1 in the method of moving surfaces, i.e. that the
``moving plane'' matrix MP is of maximal rank. Even if the second
condition fails we can still recover the implicit equation as the GCD
of the maximal minors. This is not true if the first condition
fails as illlustrated by Example 3.4.

\begin{exmp}
\label{ex2}

 Consider the surface parameterized by:

\begin{eqnarray*}
x_1 &=& 1 + s -t + st - s^2t - st^2\\
x_2 &=& 1 + s -t - st + s^2t - st^2\\
x_3 &=& 1 - s +t - st - s^2t + st^2\\
x_4 &=& 1 - s -t + st - s^2t + st^2
\end{eqnarray*}

There is one basepoint at $(s,t) = (1,1)$. The corresponding
polygon $Q$ is a pentagon with vertices $(0,0), (1,0), (0,1),
(2,1), (1,2)$. Computing the Chow form matrix gives a singular $9
\times 9$ matrix. However, we can remove one row and column to get
the $8 \times 8$ matrix below:

\begin{tiny}
\[
\begin{bmatrix}
0      & 0      & 0      & 1 -X_1 & 1 + X_1  &-1 + X_1&  1 - X_1 &-1 + X_1 \\
0      & 0      & 0      & 1 -X_2 & 1 + X_2  &-1 + X_2& -1 + X_2 &-1 - X_2 \\
0      & 0      & 0      & 1 -X_3 & -1 + X_3 &1 + X_3 & -1 - X_3 & -1 + X_3\\
1 - X_1& 1 - X_2& 1 - X_3&       0& 4X_3 - 4X_2 + 4X_1 - 4       &0       &
-4X_3 - 4X_1    & 4X_2 + 4X_3    \\
1 - X_1&-1 - X_2& -1- X_3&       0& -4 + 4X_1&0       &  8 - 4X_3 - 4X_1  &
0      \\
-1+ X_1&-1 + X_2& 1 + X_3&       0& 4 - 4X_1 &0       & 4X_3 - 4X_2 + 4X_1 - 4&- 4X_3 + 4X_1  \\
1 + X_1& 1 + X_2&-1 + X_3&       0& 0        &0       &0                  &
0       \\
-1+ X_1& 1 + X_2&-1 + X_3&       0& 0        &0       &0                  &
0       \\
\end{bmatrix}
\]
\end{tiny}

The determinant is 256 times the irreducible implicit equation which is

\begin{multline*}
2X_1-X_2+X_3-X_3^3X_1-X_2^2X_1^2+X_3X_1^2X_2-5X_3X_1+3X_2X_1\\
-2X_3X_1X_2-2X_1^2-3X_1^2X_2+2X_3X_2-2X_2^2+4X_2^2X_1+X_3X_1^2-X_2^3\\
+3X_3X_2^2-2X_3^2-X_3^2X_2+2X_2^3X_1-X_2^4-2X_3^2X_2^2-X_3X_2^3\\
+X_3^2X_1X_2+X_3^3+4X_3^2X_1
\end{multline*}

This has degree $4$ since ${\rm Area}(Q) = 5$ and there is one
basepoint of multiplicity $1$.
\end{exmp}

\begin{exmp}
\label{ex3}
Let us now consider an example where the basepoint has
multiplicity:

\begin{eqnarray*}
x_1 &=& (t + t^2)(s-1)^2 + (1+st-s^2t)(t-1)^2\\
x_2 &=& (-t -t^2)(s-1)^2 + (-1+st+s^2t)(t-1)^2\\
x_3 &=& (t-t^2)(s-1)^2 + (-1-st+s^2t)(t-1)^2\\
x_4 &=& (t+t^2)(s-1)^2 + (-1-st-s^2t)(t-1)^2
\end{eqnarray*}

Once again there is a single basepoint at $(s,t)=(1,1)$. But since,
locally the ideal at this basepoint is generated by $((s-1)^2,
(t-1)^2)$ the basepoint is an LCI. So applying the method above we
get a $15 \times 15$ matrix and an $11 \times 11$ maximal minor.

The determinant, after removing the integer constant, is

\begin{eqnarray*}
&& - 12 - 4X_1 - 9X_2 + 5X_3 - X_2^5 - 4X_3^2 - 20X_3^2X_2^3\\
&& - 16X_3^2X_2 - 32X_3^2X_2^2 - 12X_3^2X_1 - 12X_3^2X_1^2\\
&& + 8X_3X_2^4 - 12X_3^2X_1^2X_2 - 36X_3^2X_2 ^2X_1
-48X_3^2X_2X_1\\
&& + 2X_1^3 - 6X_3X_1^2 + 11X_2X_1^2 + X _1^2 - 13X_2X_1 -
3X_1^4\\
&&- 3X_1^3 X_2 + 14X_2^2X_1^2 - 9X_2^2X_1 - 16X_2^2 - 7X_1^4X_2\\
&&- X_3X _1^4 - X_1^5 - 19X_3X_1^3X_2 + 9 X_3X_1 - 15X_3X_1^3\\
&&+19X_3X_2 - 11X_2^2X_1^3 - 43X_3X_1^2 X_2 + 27X_3X_2X_1\\
&& - 6X_2^4 + 3 X_2^4X_1 + 4X_2^3X_1 - 14X_2^3 + 33X_3X_2^2\\
&&+ 10X_3X_1X_2 ^2 - 43X_3X_1^2X_2^2 + 28X_3X_ 2^3 - 12X_3X_2^3X_1
\end{eqnarray*}

The degree of this equation is $5$ and the area of the support
polygon $Q$ is 9.

\end{exmp}

\begin{exmp}
\label{ex4}
Let us now modify the above example so that the basepoints no longer form
an LCI. We will see that we can no longer recover the implicit equation
exactly from our Chow form matrix.

\begin{eqnarray*}
x_1 &=& (t + t^2)(s-1)^2 + (1+st-s^2t)(t-1)^2 + (t + st + st^2)(s-1)(t-1)\\
x_2 &=& (-t -t^2)(s-1)^2 + (-1+st+s^2t)(t-1)^2+ (t + st + st^2)(s-1)(t-1)\\
x_3 &=& (t-t^2)(s-1)^2 + (-1-st+s^2t)(t-1)^2+ (t + st + st^2)(s-1)(t-1)\\
x_4 &=& (t+t^2)(s-1)^2 + (-1-st-s^2t)(t-1)^2+ (t + st + st^2)(s-1)(t-1)
\end{eqnarray*}

Because of the additional $(s-1)(t-1)$ term, the degree of the basepoint
at $(1,1)$ drops to 3, however, the multiplicity remains $4$. Indeed,
a maximal minor of the $15 \times 15$ Chow form matrix now has rank 12.
And the determinant of {\em any} maximal minor is (up to a constant):

\begin{eqnarray*}
&&( - X_2 + 2X_3 - 1)(101 - 224X_3^5 + 8X_1 ^5 - 525X_1 + 75X_2 +
2689X_1X_3 - 573 X_3 \\
&& + 5519X_3^2X_1^2 + 3830X_3^3X _1 + 2948X_1^3X_3 + 1310X_3^2 + 155X _1X_3^2X_2^2 \\
&& - 169X_1X_3X_2^3 - 1970X_3^2 X_1^3 - 2308X_1^2X_3X_2 - 487X_3 X_2^2X_1 - 1182X_3^4X_1 \\
&& - 2296X_3^3X_1^2 + 1707X_1X_3 X_2 + 1006X_1X_3^3X_2 + 1487X_1^
2X_3^2X_2 \\
&& + 956X_1^3X_3X_2 - 1512X_3^3 - 4795X_1X_3^2 -
2118X_1X_3^2X_2 - 624X_1^4X_3 + X_2^5 \\
&& - 13X_2^4 - 88X_2^2 - 76X_2^3 - 948X_1^3X_2 + 244X_3^4X_2 - 646X_ 3^3X_2 - 513X_1X_2 \\
&& - 211X_3^2X_2^2 + 191X_2^3X_1 - 105X_1^2X_2^3 + 1140X_1^2X_2 +
185X_1^2X_3X_2^2 \\
&& + 143X_3X_2^3 + 255X_1^4X_2 + 3X_3X_2^4 - 42X_3^2X_2^3 + 19 X_1X_2^4 + 264X_1X_2^2 \\
&& - 214X_1^2X_2^2 + 48X_1^3X_2 ^2 - 385X_3X_2 + 248X_3X_2^2 + 729 X_3^2X_2 + 18X_3^3X_2^2 \\
&& + 337X_1^4 - 1050X_1^3 + 898X_3^4 - 4445X_1^2X_3 + 1133X_1^2)
\end{eqnarray*}

The second factor, of degree $5$ is the desired implicit equation.
\end{exmp}

In the last example, there is a linear extraneous factor of $-X_2
+ 2X_3 -1$. One can show that this extraneous factor divides every
maximal minor of $M_{\aa}$. Hence, the extraneous factor is
somehow intrinsic to the resultant matrix and cannot be removed.
It would be interesting to have some theoretical explanation for
this factor.

We conclude with an example where the Sylvester rows are not linearly
independent.

\begin{exmp}

\begin{eqnarray*}
x_1 &=& s + s^2 + s^3t + s^2t^2 + st^3 \\
x_2 &=& t^2(s+1) \\
x_3 &=& st(s+1) \\
x_4 &=& t(s+1) \\
\end{eqnarray*}

The Newton polygon has three interior points $st, s^2t, st^2$.  This
system turns out to have a degree $7$ LCI basepoint locus on
$X_{\aa}$.  However, one can easily check that $(s^2t - stX_3)(x_2 -
X_2x_4) = (st^2 - stX_2)(x_3 - X_3x_4)$ so that there is indeed a
syzygy of $f_1, f_2, f_3$ supported in ${\rm int}(Q)$. So there is no
maximal minor using all of the Sylvester columns. We can still
construct maximal minors using as many Sylvester columns as possible,
in this case 8 of the 9. The determinant of such a minor depends on
which choice of Sylvester columns we remove. If we remove the column
in the Sylvester block corresponding to $st \cdot f_3$ we get a matrix
whose determinant is

$$X_2(X_2^3X_1 + X_1^2X_2^2 + X_3^2X_2^2 - X_1X_2X_2 + X_2X_3^3 -
X_1X_3^2).$$

If, on the other hand, we remove a column corresponding to $s^2t \cdot f_3$
the determinant is exactly the implicit equation without the extraneous
factor of $X_2$.

\end{exmp}

\subsection{Proof of Theorem \ref{impbp}}

In this section we prove theorem \ref{impbp}. The Chow form matrix
described above, and indeed most of the formulas for Chow forms in
the literature, are applications of a general setup due to Weyman
\cite{wey}. A constructive approach using exterior algebras was
described by Eisenbud, Schreyer and Weyman \cite{ESW}. They start with an
arbitrary projective variety $X \subset \mathbb{P}^N$ of dimension
$n$ and try to compute its Chow form. Hence they consider the
incidence correspondence:

$$\begin{diagram}
    &        &   V \subset X \times G_{n+1} &      & \\
    & \ldTo^{\pi_1} &              & \rdTo^{\pi_2}& \\
    X \subset \PP^N &   &          &              & G_{n+1}
    \end{diagram}$$

Here $G_{n+1}$ is the Grassmanian of codimension $n+1$ planes in
$\mathbb{P}^N$ and $V= \{(x,F) \ : \ F(x) = 0 \}$ the incidence
subvariety of $X \times G_{n+1}$. Now given any sheaf
$\mathcal{F}$ supported on $X$ which is generically a vector
bundle, there is a complex, denoted $U_{n+1}(\mathcal{F})$ in
\cite{ESW}, of vector bundles on $G_{n+1}$ equivalent in the
derived category to $R(\pi_2)_{\ast}\pi_1^{\ast}\mathcal{F}$. This
leads to the following completely general result.

\begin{thm} \label{general} Let $X \subset \mathbb{P}^N$ be any variety of
dimension $n$. Let $\mathcal{F}$ be any sheaf supported on $X$
that is generically of rank 1. Let $F_0, \dots, F_n$ be any
linearly independent sections of $\mathbb{P}^N$ which
simultaneously meet $X$ only at finitely many points at all of
which $\mathcal{F}$ is of rank 1. The last map in the complex
$U_{n+1}(\mathcal{F})$ has cokernel of rank equal to the degree of
the zero-dimensional subscheme of $X$ cut out by $F$.
\end{thm}

\begin{pf}

Consider again the incidence correspondence.  As
$U_{n+1}(\mathcal{F})$ is isomorphic in the derived category to
$R(\pi_2)_{\ast}\pi_1^{\ast}\mathcal{F}$, the cokernel of the last
map in particular is just $(\pi_2)_{\ast}\pi_1^{\ast}\mathcal{F}$
itself. So all we need to show is that the dimension of the the
fiber of this sheaf at a point $F \in G_{n+1}$ satisfying the
above properties is the degree of the subscheme $X_F$ of $X$
defined by $F$.

First consider the fiber of the morphism $\pi_2$ over $F$. Let $R$
be the coordinate ring of $X$ and $S$ the Stiefel coordinate ring
of $G_{n+1}$ with variables $\underline{a}$. The ideal of $V$ in
$R \otimes S$ is denoted $I(\underline{a})$. Now, by definition
the fiber over the point $F$ defined by a choice $\underline{a} =
a$ with corresponding maximal ideal $m_a$ in $S$ is $(R \otimes
S)/I(\underline a) \otimes_{S} S/m_a$. But this is just $R/I(a)$
which is the coordinate ring of $X_F$. Hence the fiber of $\pi_2$
over $F$ is $X_F \times F$. (Note that different choices of $a$
realizing the same point $F$ give the same ideal $I(a)$).

Next, since $X_F$ is a zero dimensional subscheme of the generic
locus of $\mathcal{F}$ it is actually affine and $\mathcal{F}$ is
trivial on $X_F$. Let $R/I(a)$, as above, be the (dehomogenized)
coordinate ring of $X_F$ and hence also of $X_F \times F$. As our
sheaf was trivial, the pushforward onto the closed point $F$ is
just $R/I(a)$ itself viewed as a vector space over the residue
field of $F$. The dimension of this vector space is by definition
the degree of $X_F$ as desired.

\end{pf}

We can now prove Theorem \ref{impbp} as a corollary.

\begin{pf}

We consider, in this case, $\mathcal{F} = \oo({\rm int}(2Q))$ the
divisor corresponding to the interior of the polytope $2Q$. In
\cite{khe}, it was shown that $U_{3}(\mathcal{F})$ reduced to a two
term complex with matrix exactly as described above. The sheaf
$\mathcal{F}$ is locally free of rank 1, except possibly on the
singular points of $X_\aa$, which only occur on the vertices of
$Q$. By the construction of $\aa$, at least one of $x_1, \dots, x_4$ does
not vanish on each vertex, hence the base point locus always misses
the singular locus.

Now, we can apply Theorem \ref{general}. Pick a maximal minor of
our matrix. For a generic $X_1, X_2, X_3$ not on the surface $S$,
this remains a maximal minor of the specialization. Moreover, the
corank of this minor is the degree of $I(f_1(X_1), f_2(X_2),
f_3(X_3))$. However, for a point $X_1, X_2, X_3$ on the surface,
the number of basepoints increases, therefore the rank of our
matrix $M$ decreases, hence the determinant of our chosen minor
must be zero. Moreover, the rank drop of the minor for a generic
point on the surface is exactly $\deg{\phi}$ (the number of ``new
basepoints'' mapping on to our point).  Since any order $k$
derivative of the determinant of a matrix of linear forms is in
the ideal of corank $k$ minors (easy to see from the expansion of
determinant), the first $\deg{\phi}-1$ derivatives of the
determinant are also zero for a generic point on the surface.
Since $P$ was irreducible, $P^{\deg{\phi}}$ must divide our chosen
maximal minor.

For the second part, in the case of an LCI, the corank of our maximal
minor, i.e. the degree of the base point locus, is the same as the sum
of the multiplicities of our base points. If moreover, the maximal
minor is chosen to contain all Sylvester rows and columns, only
B\'ezout rows and columns are removed, each of which drops the degree
by 1.  Thus the degree of our determinant is equal to the degree of
$P^{\deg{\phi}}$ and so they must be equal up to a constant.
\end{pf}

\section{The method of moving surfaces}\label{s:4}

We now switch gears and present an entirely different method for
constructing matrix formulas in implicitization. For the rest of
this paper we will work with the projective surface $\Phi \subset
\PP^3$ defined by the four coordinates $X_1, X_2, X_3, X_4$.

The idea will be to construct linear and quadratic syzygies on the
polynomials $x_1, x_2, x_3, x_4$ and put them together into a matrix
of linear and quadratic forms in the $X_i$. For the case of
homogeneous and bihomogeneous polynomials, this is the method of
moving planes and surfaces introduced by Sederberg and Chen \cite{SC}.
However the proof we present in Section 5 is quite different, and in
our opinion more insightful, than the ones in the literature. Our goal
will to be to extend the method to general toric surfaces which will
require looking at certain ``degrees'' of the homogeneous coordinate
ring of the toric variety.

We shall see that the syzygy method has certain advantages and
disadvantages to the Chow form/resultant method described above.  It
will always give smaller matrices due to the fact that some of the
entries are quadratic in the $X_i$. Second, the algorithm will be
relatively easy to describe and efficient in practice; all of the
computations are just numerical linear algebra. Finally, the method
appears to be surprisingly flexible in the presence of base points. We
shall see empirical evidence supporting this at the end of the
section.

On the other hand, rigorous proofs of the method in any of the more
complicated situations have been hard to come by. Also, as pointed out
above, all of the computations are linear algebra in the coefficients
of the $x_i$. In particular, the method becomes much more inefficient
with a generic parameterization or whenever the coefficients of the
$x_i$ are not numerical. The Chow form matrix constructed above, on
the other hand, works the same for arbitrary coefficients and is
therefore preferred when implicitizing a family of surfaces.

\subsection{Moving planes and quadrics}

Given a rational surface $\Phi$ parameterized by

\begin{equation}\label{par2}
\left\{
\begin{array}{ccc}
X_1&=&x_1(s,t) \\
X_2&=&x_2(s,t) \\
X_3&=&x_3(s,t) \\
X_4&=&x_4(s,t)
\end{array}
\right.
\end{equation}

A {\em moving plane} is a syzygy on $I = \langle x_1, x_2, x_3, x_4 \rangle$,
i.e an equation of the form

$$A_1(s,t)X_1 + A_2(s,t)X_2 + A_3(s,t)X_3 + A_4(s,t)X_4$$

which is identically zero as a polynomial in $s$ and $t$ after the
specialization $X_i\mapsto x_i.$ Notice that each particular choice of
$(s,t)$ gives the equation of a plane which intersects the surface
$\Phi$ at the point $(x_1(s,t), \dots, x_4(s,t))$. Hence, this is said
to be a plane that follows the surface $\Phi$ and justifies the
terminology moving plane.

Similarly, a {\em moving quadric} is a syzygy on $I^2$:

$$A(s,t)X_1^2 + B(s,t)X_1X_2 + \cdots + J(s,t)X_4^2$$

Once again a choice of $(s,t)$ gives the equation of a quadric
meeting the surface $\Phi$. Hence, the moving quadric is said to
follow the surface.

If we rewrite the moving planes and quadrics in terms of the
monomial bases in $s$ and $t$ we get vectors of linear or
quadratic forms in the $X_i$. Clearly multiplying each moving
plane by $X_1, X_2, X_3, X_4$ gives a moving quadric. Therefore,
we will only look for ``new'' moving quadrics. If we can now get
enough of these vectors, we may be able to build a square matrix
out of them. The determinant of this square matrix will hopefully
be equal to the implicit equation of $S$. The following well known
result is our starting point.

\begin{prop} \label{SCthm}
Let $M(X_1, X_2, X_3, X_4)$ be any square matrix constructed from
moving planes and quadrics as above. Then $\det (M(x_1, x_2, x_3,
x_4)) = 0$. In particular the implicit equation always divides the
determinant of $M$ (which may, quite possibly, be identically 0).
\end{prop}

\begin{pf}
This has been proved in even more generality in \cite{SC}.
\end{pf}

The big question is now, of course, how should the moving planes
and quadrics be chosen? In the case of homogeneous polynomials they
were chosen to also be homogenous of an appropriate degree. In the
case of bihomogeneous polynomials, the moving planes and quadrics
can  be chosen to be bihomogeneous. In the more general toric
setting we will need to work in appropriate homogeneous
coordinates for the set $\aa$.

\subsection{Homogeneous coordinate ring of $X_{\aa}$}

Let $\aa$ be the union of monomials in the $x_i$ as before and $Q
= {\rm conv}(\aa)$ the associated polygon. Let $E_1, \dots, E_s$
be the edges of $Q$ and $\eta_1, \dots, \eta_s$ the corresponding
inner normals. 

We can therefore define $Q$ by its facet inequalities.

$$Q = \{ m \in \mathbb{R}^2 \ : \ \langle m, \eta_i \rangle \geq -a_i \ {\rm for} \ i= 1, \dots s \}$$

\noindent For some $(a_1, \dots, a_s) \in \mathbb{Z}^s$.

$X_{\aa}$ is a toric variety with a given very ample line bundle
determined by the polytope $Q$. We will need to consider other
divisors on $X_{\aa}$. David Cox \cite{cox} defined a single ring that
encapsulates all torus invariant divisors on $X_{\aa}$. 

\begin{defn}
The {\em homogeneous coordinate ring} for $X = X_{\aa}$ is the
polynomial ring $S_X = \KK[y_1, \dots, y_s]$ where the monomials are
graded as described below.
\end{defn}

Consider the exact sequence of maps:

$$ 0 \to \mathbb{Z}^2 \overset{\phi} \to  \mathbb{Z}^s  \overset{\pi}{\to}  G \to
0$$

Here $\phi$ is the map $m \to (\langle m, \eta_1 \rangle, \dots,
  \langle m, \eta_s \rangle) $.  The ring $S_X$ is graded by
  elements of $G$ where $\deg{y^{\alpha}} = \pi(\alpha)$.

The graded pieces of this ring have bases corresponding to lattice
points in polygons.  More precisely the monomials in $S_{\pi(b)}$ are
in one to one correspondence with the lattice points in $Q_b = \{ m
\in \mathbb{R}^2 \ : \ \langle m, \eta_i \rangle \geq -b_i \}$. And
moreover, $\pi(b) = \pi(b')$ iff $Q_b$ is a translate of $Q_{b'}$

So it will make sense to talk about $S_{Q_b}$, the graded piece of $S$
defined by $Q_b$.

\begin{rem}
In truth the divisors and homogeneous coordinate ring are really
defined for the {\em normal} toric variety $X_Q$ obtained from the
normal fan of $Q$.  This variety is the normalization of our $X_{\aa}$.
The projection and all prior and subsequent results can be lifted up
to $X_Q$ without affecting any of the calculations.
\end{rem}

\subsection{Picking moving planes and quadrics}

Also associated to the polygon $Q$ is a certain polynomial $E(k)$, the
Ehrhart polynomial defined in \cite{sta}, which counts the number of
lattice points in $k \cdot Q$.  In the case $Q$ is two dimensional, it
turns out that

$$E(x)=Ax^2+\frac{B}{2}x+1$$

\noindent where $A=\frac{Area(Q)}{2}$ and $B$ equals the number of
boundary points.

Let $I$ be a nonempty proper subset of $\{1, \dots, s\}$ such that the
corresponding edges form a connected set. Let $E_I$ be this connected
set of edges of $Q$, let $B_I$ be the sum of the lattice edge lengths
of $E_I.$ It is easy to see that the number of lattice points in the
set of edges $E_I$ in $k \cdot Q$ is $B_I k +1$.

\begin{assum} \label{ass1}
We choose $E_I$ in such a way that $B\geq 2B_I.$
\end{assum}
\begin{rem}
Observe that this can always be done, for instance by taking as $E_I$
the shortest edge of $Q.$ In practice, we will want to pick $E_I$ in
such a way that $B_I$ is as big as possible consistent with Assumption
\ref{ass1}.
\end{rem}

Now we can define a degree of $S$ denoted $S_{Q\setminus E_I}$
obtained by ``pushing in'' all 
of the edges of $Q$ in $E_I$ by one,
whose monomial basis consists of all lattice points in $Q$ not on any
of the edges $E_I$.  In the case of homogenous polynomials of degree
$n$, the only $E_I$ satisfying Assumption \ref{ass1} consist of a
single edge and the degree in question in just $n-1$. In the case of
bihomogeneous polynomials of bidegree $(m,n)$, we can take $E_I$ to be
two consecutive edges and the degree is $(m-1, n-1)$. Note that in the
latter case $B - 2B_I = 0$ which, as we shall see, means that we will
not need to take any moving planes and can build a matrix entirely out
of moving quadrics. We now formally define what we mean by moving
planes and quadrics of this degree.

Consider the following $\KK$-linear map
\begin{equation}\label{MP}
\begin{array}{cccc}
\psi_1:&{S_{Q\setminus E_I}}^4&\to&S_{2Q\setminus E_I} \\
&(p_1,p_2,p_3,p_4)&\mapsto&\sum_{i=1}^4 p_i x_i,
\end{array}
\end{equation}
and let $MP$ be the matrix of this map in the monomial bases.

\begin{defn}
As in \cite{CGZ}, any element of the form
$(A_1,A_2,A_3,A_4)\in\ker(\psi)$ will be called a \textit{moving
plane} of ``degree'' $Q\setminus E_I$ that follows the surface
\eqref{par}. Sometimes, we will write moving planes as
$A_1X_1+A_2X_2+A_3X_3+A_4X_4.$
\end{defn}

Now for moving quadrics we consider the following map:

\begin{equation}\label{mq}
\begin{array}{cccl}
\psi_2:&{S_{Q\setminus E_I}}^{10}&\to&S_{3Q\setminus E_I} \\
&(A_{i,j,k,l})_{i+j+k+l=2}&\mapsto&\sum_{i+j+k+l=2}A_{i,j,k,l}x_1^ix_2^jx_3^kx_4^l,
\end{array}
\end{equation}
and let $MQ$ be the matrix of $\psi_2$ in the monomial bases. Then

$$\begin{array}{l}
\#\,\mbox{rows of}\,MQ = \#\left(3Q\setminus E_I\right)\cap\ZZ^2=\\
=(9A+\frac{3}{2}B+1)-(3B_I+1)=9A+\frac{3}{2}B-3B_I
\end{array}
$$
and
$$\begin{array}{l}
\#\,\mbox{columns of}\, MQ = 10\,\#\left(Q\setminus E_I\right)\cap\ZZ^2 =\\
=10\left(A+(B-2B_I)/2\right)=10A+5B-10B_I.
\end{array}$$

Now a {\em moving quadric} of degree $Q \setminus E_I$ which
follows our surface $S$ is just an element of the kernel of $MQ$.

We now describe the method of moving quadrics. It differs from the
presentations in the literature not only in its application to general
toric surfaces but also in that we allow the bases of moving planes
and quadrics to be chosen freely. Earlier papers specify that moving
planes and quadrics be chosen of a specific form to ensure that the
resulting matrix has determinant non-zero. Our more intrinsic proof
of Section 5.1 makes this unnecessary.

\begin{itemize}

\item Compute a basis $P$ of the kernel of $MP$. The entries are $P_i
= A^i_1 X_1 + A^i_2 X_2 + A^i_3 X_3 + A^i_4 X_4$. Where the $A^i_j$ are
polynomials in $S_{Q\setminus E_I}$.

\item Each $P_i \cdot X_j$ for $j = 1, \dots, 4$ is in the kernel of
$MQ$. We will see that these are linearly independent. Extend this set
to an entire basis for the kernel of $MQ$.  Let $Q_1, \dots, Q_d$ be
the new moving quadrics in this basis.

\item Construct a matrix $\MM$ out of the $P_i$ and $Q_j$ such that the
columns correspond to the monomial basis of $S_{Q\setminus E_I}$ and
the entries are the linear (or quadratic) polynomial in $X_1, \dots,
X_4$ corresponding to the coefficient of that monomial in $P_i$ (or
$Q_j$).

\end{itemize}

Our hope is that the resulting matrix will be square and that the
determinant is the implicit equation. To start with, by Theorem
\ref{SCthm}, if the matrix $\MM$ has more rows than columns, then
the determinant of any maximal minor (possibly 0) is divisible by
the implicit equation.

\section{Validity of the method of moving quadrics without basepoints}

In this section we verify, in the absence of basepoints, that the
method of moving quadrics gives a square, nonsingular matrix whose
determinant is exactly the implicit equation raised to the power
the degree of the parameterization. We will need to make one
assumption:

\begin{assum}
The moving plane matrix $MP$, or the map $\psi_1$, has maximal rank.
\end{assum}

This assumption also appears in the papers by Cox, Goldman, and
Zhang \cite{CGZ} and D'Andrea \cite{dan}. Empirical evidence
suggests that it is almost always satisfied.  It appears that for
a fixed $Q$ and $E_I$, and any generic set of $x_1, x_2, x_3, x_4$
without basepoints, $MP$ has maximal rank.

 We now build a complex containing both the moving plane map
$\psi_1$ and the moving quadric map $\psi_2$ in Figure \ref{complex}.

\begin{figure} 
$$\begin{diagram}
   &      &        &      &  0            &                   & 0             &          &   \\
   &      &        &      &  \dTo         &                   & \dTo          &          &   \\
   &      & 0      & \rTo & (S_{\QEI})^6  &  =                & (S_{\QEI})^6  &     \rTo & 0  \\
   &      & \dDashto &    &  \dTo^i       &                   & \dTo^{x'}     &   &          &    \\
 0 & \rTo & K_1^4  & \rTo &(S_{\QEI})^{16}& \rTo^{{\psi_1}^4} & (S_{2\QEI})^4 & \rTo     & 0  \\
   &      & \dTo^X &      &  \dTo^{X}     &                   & \dTo^x        &          &    \\
 0 & \rTo & K_2    & \rTo &(S_{\QEI})^{10}& \rTo^{{\psi_2}}   & S_{3\QEI}     & \rDashto & 0  \\
   &      & \dTo   &      &  \dTo         &                   & \dTo          &          &    \\
   &      & \tilde{K_2}&   &   0           &                   &  0            &          &
\end{diagram}$$
\caption{Complex of moving planes and quadrics}
\label{complex}
\end{figure}

The terms $K_1$ and $K_2$ are the kernels of the moving plane and
moving quadric maps $\psi_1$ and $\psi_2$ respectively. The term
$\tilde{K_2}$ is the cokernel of the map $X$ of $K_1^4$ into
$K_2$, generated precisely by a basis of $K_2$ extending the image
of moving planes multiplied by linear forms. In this new language
a matrix $\MM$ of moving planes and quadrics is a basis for $K_1^4
\oplus \tilde{K_2}$ taken as vectors in $(S_{\QEI})$ with
coefficients that are linear or quadratic forms in $(X_1, X_2,
X_3, X_4)$.

We now prove out two main theorems that together prove the
validity of the method of moving quadrics.

\begin{thm}\label{square}
If $MP$ has maximal rank, then $\dim(K_1) + \dim(\tilde{K_2}) =
\dim(S_{\QEI})$ and $\dim(K_1) + 2\dim(\tilde{K_2}) =
\rm{Area}(Q)$. Therefore, the method of moving quadrics yields a
square matrix with determinant of degree equal to the implicit
equation.
\end{thm}

\begin{thm}\label{nonsing}
Let $p = (p_1,p_2,p_3,p_4)$ be a point not on the surface $X$. The
moving plane matrix $\MM$ is nonsingular at $p$. Consequently, if
$\psi_1$ has maximal rank $\det(\MM) = P^{\deg(\phi)}$ where $P$
is the implicit equation as desired.
\end{thm}

Before proceeding we further describe the maps in the complex. The
second row consists of four copies of the moving plane complex. An
element of $(S_{\QEI})^{16}$ is represented as a four tuple of
linear forms in $X_1, X_2, X_3, X_4$ with coefficients in
$S_{\QEI}$. Similarly the bottom row is the moving quadric
complex. An element of $(S_{\QEI})^{10}$ is a quadratic form in
$X_1, X_2, X_3, X_4$  with coefficients in $S_{\QEI}$. generated
by the 10 monomials $X_iX_j$ with $i \leq j$ The map $X$,
multiplication by $(X_1, X_2, X_3, X_4)$, sends the four tuple
$(u_1, u_2, u_3, u_4)$ of linear forms to the quadratic form $\sum
u_i X_i$. This has the effect of sending $X_i$ in position $j$ and
$X_j$ in position $i$ both to $X_iX_j$.

The kernel of $X$ is isomorphic to $(S_{\QEI})^6$ indexed by pairs
$(i,j)$ with $i < j$. The injection $i$ sends the term $p_{ij}$ to
$(0,\dots, p_{ij}X_j \dots, -p_{ij}X_i, \dots 0)$ with  $X_j$ in
position $i$ and $-X_i$ in position $j$. The rightmost column is a
graded piece of the Koszul complex on $(x_1, x_2, x_3, x_4)$, with
$x$ mapping a four tuple $(s_1, s_2, s_3, s_4)$ to $\sum s_ix_i$
and $x'$ sending $p_{ij}$ with $i < j$ to $(0, \dots, p_{ij}x_j,
\dots, -p_{ij}x_i \dots, 0)$.

Commutativity of the diagram is immediate. The rows are all exact
by construction. The second column is also clearly exact. The
rightmost column is more interesting. When $(x_1, x_2, x_3, x_4)$
have no basepoints, the map $x'$ is injective and $x$ is
surjective. This can be seen by investigating the complex
$U_4(\oo(3\QEI)))$ arising from the Tate resolution in the theory
of \cite{ESW}. However, the spot in the middle is not exact. We
shall see later that obstruction to exactness comes from a certain
'Bezoutian' map determined exactly by elements of $\tilde{K_2}$.

Now, to prove Theorem \ref{square} we will need three lemmas:

\begin{lem}\label{l1}
If $MP$ has maximal rank, then the number of linearly independent
moving planes of degree $Q\setminus E_I$ which follow the surface
is $B-2B_I.$
\end{lem}

\begin{pf}
There are $$(4A+B+1)-(2B_I+1)=4A+B-2B_I$$ integer points in
$2Q\setminus E_I,$ and
$$(A+\frac{B}{2}+1)-(B_I+1)=A+(B-2B_I)/2$$ integer points in
$Q\setminus E_I.$ If $MP$ has maximal rank, then the number we
want to compute is the dimension of the kernel of $\psi_1$ which
equals
$$4\left(A+(B-2B_I)/2\right)-(4A+B-2B_I)=B-2B_I$$
as claimed.
\end{pf}

\begin{lem} \label{l2}
If $\psi_1$ is surjective then so is $\psi_2$.
\end{lem}

\begin{pf} This is an easy diagram chase. Given $s \in S_{3\QEI}$  pull
it back to $S_{2\QEI}^4$ and then to $(S_{\QEI})^{16}$ via the
surjectivity of the corresponding maps. Finally map this down to
$t \in (S_{\QEI})^{10}$. Commutativity of the diagram yields
$\psi_2(t) = s$.
\end{pf}

\begin{lem} \label{l3}
The map $X$ from $K_1^4$ to $K_2$ is injective.
\end{lem}

\begin{pf} Given $k$ in the kernel, it is a non-zero
element of $(S_{\QEI})^{16}$ mapping to zero in $(S_{\QEI})^{10}$.
By exactness it comes from a nonzero element in $S_{\QEI}^6$
mapping to a nonzero element $k'$ in $(S_{2\QEI})^4$. But
commutativity implies $k' = \psi_1(k) = 0$, a contradiction.
\end{pf}

We are now ready for the proof of  Theorem \ref{square}.

\begin{pf}
By Lemma \ref{l2} $MQ$ has maximal rank. From the computations of
the last section, the dimension of $(S_{\QEI})^{10}$ is
$10A+5B-10B_I$, while the dimension of $S_{3\QEI}$ is $9A +
\frac{3}{2}B - 3B_I$. Thus the rank of $K_2$ is $A + \frac{7}{2}B
- 7B_I$. By Lemma \ref{l1} the rank of $K_1^4$ is $4(B-2B_I)$, so
by Lemma \ref{l3} the rank of $\tilde{K_2}$ is $A - \frac{B}{2} +
B_I$.

So, the sum of the ranks of $K_1$ and $\tilde{K_2}$ is

\begin{eqnarray*}
 B - 2B_I + (A - \frac{B}{2} + B_I) &= A + \frac{B}{2} - B_I \\
                                    &= \dim S_{\QEI}
\end{eqnarray*}

Moreover, the total degree of the determinant is

$$ B - 2B_I + 2(A - \frac{B}{2} + B_I) = 2A $$

This is twice the Euclidean area, hence equal to the normalized
area of $Q$ as desired.
\end{pf}

The theorem just proved shows that $\MM$ is square of the right
rank. Theorem \ref{nonsing} will show that its determinant does
not vanish outside of the surface.

Let  $p = (p_1, p_2, p_3, p_4) \in \PP^3$ be a point not on the
surface parametrized by $X$. WLOG assume that $p_4 = 1$. Make a
change of coordinates $X_1' = X_1 - p_1X_4, X_2' = X_2 - p_2X_4,
X_3' = X_3 - p_3X_4$ and $X_4' = X_4$. The point $(p_1, p_2, p_3,
p_4)$ is transformed to $(0,0,0,1)$. Since the parameterization
has no base points, the $Ch_\aa(x_1', x_2', x_3') \neq 0$ by
Theorem \ref{chow}.

We now use two facts arising from resultant complexes.

\begin{lem} \label{l4}
The restricted map $\tilde{\psi_1} \ : \ S_{\QEI}^3 \to S_{2\QEI}$
given by $(s_1, s_2, s_3) \to \sum s_i x_i'$ is injective. In particular
no moving plane $A_1X_1' + A_2X_2' +A_3X_3' +A_4X_4'$ vanishes
at $(X_1',X_2', X_3', X_4') = (0,0,0,1)$.
\end{lem}

\begin{pf}
In \cite[Theorem 3.4.1]{khe}, a matrix whose determinant gives
$Ch_\aa(x_1,x_2,x_3)$ is constructed, and this matrix has a Sylvester
part coming from $\tilde{\psi_1}$. As we have
$Ch_\aa(x_1',x_2',x_3')\neq0,$ it turns out that $\tilde{\psi_1}$
must be injective. Any vanishing moving plane as above has $A_4 = 0$
so must in fact be in the kernel of $\tilde{\psi_1}$.
\end{pf}

\begin{lem} \label{l5}
The restriction of the last column:

$$ 0 \to S_{\QEI}^3 \to S_{2\QEI}^3 \to S_{3\QEI} \to, 0$$

which is just the Koszul complex on $x_1', x_2', x_3'$, is exact.
\end{lem}

\begin{pf}
 We may consider $x_1',x_2',x_3'$ as sections of sheaves on the toric
variety $X_\AA.$ As in \cite[Section 4]{DK}, we start with the
Koszul complex of these sheaves in degree $3\beta-\beta_I,$ where
$\beta$ is the degree associated to $Q\cap\ZZ^2$ and $\beta_I$ the
divisor associated to all the edges whose union equals $E_I.$ As
in the proof of Theorem $3.1$ in \cite{DK}, one can see that we
can apply the Weyman's complex (see \cite[Section 3.4.E]{GKZ}) to
this complex. By the toric version of Kodaira vanishing (see
\cite{mus}), all higher cohomology terms vanish and we get that
the complex above is generically exact. Indeed, the determinant of
the complex equals the Chow form of $x_1',x_2',x_3'.$
\end{pf}

\begin{cor}\label{c1}
Any moving quadric $\sum_{1 \leq i \leq j \leq 4} A_{ij} X_i'X_j'$
vanishing at $(0,0,0,1)$ is in the image of $K_1^4$ under the
multiplication map $X$.
\end{cor}

\begin{pf}
Start with such a vanishing moving quadric $q$. Plugging in we see
that $A_{44} = 0$. Hence $q = q_1X_1' + q_2X_2' + q_3X_3'$ where $q_1
= A_{11}X_1' + A_{12}X_2' + A_{13}X_3' + A_{14}X_4'$, $q_2 =A_{22}X_2'
+ A_{23}X_3' + A_{24}X_4$ and $q_3 = A_{33}X_3' + A_{34}X_4'$.

Pulling back to $q' = (q_1, q_2, q_3, 0) \in (S_{\QEI})^{16}$ and
mapping to $(S_{2\QEI})^4$ by substituting $x_1', x_2', x_3'$ into
$q_1, q_2, q_3$ we get an element of the subspace $(S_{2\QEI})^3$
as in the restricted complex above which is still in the kernel of
$X$.  Thus, by Lemma \ref{l2} we can pull back via $X'$ to
$S_{\QEI}^3 \subset S_{\QEI}^6$. Let $q''$ be the image of this
element in $S_{\QEI}^{16}$. By construction $\psi_1^4(q' - q'') =
0$ and $X(q' - q'') = q$. But now we can pull back $q' - q''$ to
$k = (k_1, k_2, k_3, k_4)$ with $X(k) = q$ as desired.
\end{pf}

It is now straightforward to finish the proof of Theorem
\ref{nonsing}.

\begin{pf} Suppose $(u(p),v(p))$ is in the kernel. Write
$u = \sum u_i X_i'$ and $v = \sum v_{ij} X_i'X_j'$.  Substituting
in for $p$ we have $u_4 + v_{44} = 0$.  Therefore the moving
quadric $X_4'u + v$ has no $(X_4')^2$ term and thus vanishes at
$p$. By Corollary \ref{c1} this must be in the image of $K_1^4$ so
we must have $v = 0$. But now $u(p) = 0$ violating Lemma \ref{l4}.
Hence $\MM$ is singular only on points of $X$.  If $\psi_1$ is
maximal rank then $\MM$ is square, hence its determinant is a
power of the implicit equation. Since the degree of $\det(\MM) =
Area(Q)$, the exponent must be $\deg(\phi)$.
\end{pf}

\begin{exmp}

Consider the system from Example \ref{ex1}:
\begin{eqnarray*}
x_1 &=& s^3 + t^2\\
x_2 &=& s^2 + t^3\\
x_3 &=& s^2t + st^2\\
x_4 &=& st
\end{eqnarray*}

The total boundary length of the quadrilateral $Q$ is 7. We can
pick $E_I$ to be the long edge of length 3. Hence $B - 2B_I = 1$.
Applying the method of moving quadrics then gives a matrix with
one moving plane and two moving quadrics:

$$ \left[ \begin{array}{ccc}
 -X_1 - X_2 - X_3 & X_3 + X_4 & X_3 + X_4 \\
&& \\
X_1X_3 - X_2X_4 + X_4^2 & X_1X_3 - X_2X_4 - X_3X_4 & -X_3^2 + X_4^2 \\
&& \\
-X_1^2 - X_1X_2 - 3X_1X_3 &  X_1X_4 + X_2X_4 + X_3X_4 & X_1X_2 + X_2X_3 + 2X_3^2 \\
+ 2X_2X_4 - X_3^2 + X_3X_4 & & -X_3X_4 - 2X_4^2
\end{array} \right]
$$

The determinant is exactly the degree $5$ implicit equation.

\end{exmp}

\subsection{Moving quadrics in the presence of base points}

In the case of homogeneous parameterizations ($X_{\aa} = \PP^2$),
\cite{BCD} gives a series of conditions for when the method of moving
quadrics works even with basepoints. The conditions are labelled
(BP1)-(BP5) but essentially they boil down to assuming the basepoints
form an LCI, there are no syzygies on linear combinations of $x_1,
x_2, x_3, x_4$ of the desired degree, and that there are the ``right
number'' of moving planes of the degree in question.

The last assumption can be rephrased into a regularity assumption on
the ideal of basepoints $I$. Using commutative algebra on graded rings
they deduce a corresponding regularity bound on $I^2$ which implies
that there are also the ``right number'' of linearly independent
moving quadrics.

To extend these conditions to the toric setting would seem to
require a notion of ``toric regularity'' using the homogeneous
coordinate ring $S_X$ in place of the usual graded polynomial
ring. Perhaps the definition proposed by Maclagan and Smith
\cite{MS1,MS2} can be applied here. Instead of delving into the
theory of toric commutative algebra and what does and does not
extend, we simply present some examples to illustrate how the
toric method of moving quadrics can often work in the presence of
basepoints.

\begin{exmp} We repeat Example \ref{ex2} using moving quadrics.

\begin{eqnarray*}
x_1 &=& 1 + s -t + st - s^2t - st^2\\
x_2 &=& 1 + s -t - st + s^2t - st^2\\
x_3 &=& 1 - s +t - st - s^2t + st^2\\
x_4 &=& 1 - s -t + st - s^2t + st^2
\end{eqnarray*}

Recall that we have one basepoint at $(1,1)$ with multiplicity 1.  If
there were no basepoints then we can choose $B - 2B_I = 1$ and we
would expect one moving plane and two moving quadrics. Applying
the algorithm gives two planes and one quadric but still a $3 \times 3$
square matrix:

$$\begin{bmatrix}
-X_3 + X_4 & 0 & X_1 - X_2 \\
X_2 - X_3 + 2X_4 & X_2 + X_3 & -X_2 - X_3 + 2X_4 \\
X_2X_1 + X_3X_1 & X_3X_1 - X_1X_4 + X_2^2 + X_2X_4 & -2X_1^2 + X_2^2 + X_2X_4 - X_3X_4 + X_4^2
\end{bmatrix}$$

The method of moving quadrics works perfectly here and gives the implicit
equation of degree 4.
\end{exmp}

Example \ref{ex3} which had an LCI basepoint of multiplicity 4 also
works with the method of moving quadrics. In this case we get
$5$ moving planes and no moving quadrics. The implicit equation is
recovered as the determinant of the corresponding $5 \times 5$ matrix
of linear forms.

Example \ref{ex4} has a basepoint which is not an LCI. In this case,
the moving quadric matrix was not square. Indeed there were four moving
planes and two moving quadrics on a space of five monomials.

However, taking the maximal minor consisting of the four planes
and either one of the two quadrics gives the implicit equation
with a linear extraneous factor. Unlike, the Chow form matrix of
Example \ref{ex4}, this extranous factor is not intrinsic to the
construction. The two different maximal minors give different
extraneous factors, hence the implicit equation is the gcd of the
maximal minors.

\section{Conclusion}

In this paper we extend two of the most important implicitization
techniques, resultants and syzygies, to general toric surfaces.
There are a couple of interesting open questions remaining.

For the resultant method, when the basepoints are not an LCI
every maximal minor of the resultant matrix will have an extraneous
factor. Is there a way to compute this extraneous factor apriori?

For the syzygy method, the biggest open question is how to
extend the method when basepoints are present. Our examples
show that the method may often still work. The second
open problem is an understanding of exactly when the
moving plane matrix has maximal rank.

\section*{Acknowledgements}

We would like to thank Bernd Sturmfels for directing us
to this problem and David Cox for going over earlier manuscripts
with a fine tooth comb.

\end{document}